%% file: SurveyLaplacian.tex
\documentclass[12 pt,amstex]{article}

\usepackage[usenames]{color}
\usepackage{url}
\usepackage{pgf,tikz}
\usetikzlibrary{arrows}
\usepackage{amsmath}
\usepackage{amssymb}

\usepackage{authblk}

\input{inputlatex}

\begin{document}



\title{A Primer on Laplacian Dynamics in Directed Graphs}

\input{authorInfo-arXiv}

\date{\today}


\maketitle



\section{Introduction}
\label{chap:intro}

Directed graphs (or digraphs) are an important generalization of undirected graphs and they have wide-ranging
applications.  Examples include models of the internet \cite{Brod} and social networks
\cite{Carr},  food webs \cite{May}, epidemics \cite{Jom}, chemical reaction networks \cite{Rao},
databases \cite{Ang}, communication networks \cite{Ahls}, the Pagerank algorithm \cite{Stern},
and networks of autonomous agents in control theory \cite{Fax} to name but a few.
In many of these applications, it is of crucial importance to understand the asymptotics
(as $t\rightarrow \infty$) of solutions of the first order Laplacian differential equation
\bse
\dot x=-Lx \quad \logand \quad \dot p = -pL .
\label{eq:diffeq}
\ese
We will show how the first these equations is associated with the physical process of consensus
and the second with diffusion. We will give a unified treatment of both of these in which
one is treated as the dual to the other. We also show how this extends to the discrete
versions of these processes.

We describe the basic theory of Laplacian dynamics on directed graphs that are weakly connected.
The restriction of this theory to \emph{undirected} graphs is well documented in
textbooks (see \cite{Godsil2001_AlgGraphTheory}, \cite{Chung1997_SpectralGraphTheory}),
but as far as we know, this is the first complete exposition of the
general theory (\emph{directed} graphs) in  a single work.

Many of the results we will discuss had earlier been ``folklore" results living largely outside
the mathematics community and not always with complete proofs (see \cite{caughveer,veerkumm}
for some references). In the mathematics community, directed graphs are still much less studied
than undirected graphs (especially true for the algebraic aspects). As a consequence, there are
not many good mathematics books on the subject.

Part of the reason for that is probably that directed graphs are a lot messier than undirected graphs.
For example, we will see that while \emph{undirected} graphs are either connected or not,
for \emph{directed} graphs) there are various gradations of connectedness.
Another complication is that while Laplacians of \emph{undirected} graphs are diagonalizable
and have real eigenvalues, neither statement is necessarily true for Laplacians of
\emph{digraphs}. Thus, many statements for undirected graphs take more work to prove, or are wrong.

Another reason for confusion is that there is no standard way to orient a graph.
The in-degree Laplacian of $G$ is the same as the out-degree Laplacian for $G'$, the graph $G$
with all orientations reversed.  In \cite{veerkumm}, the convention was proposed where the direction
of edges corresponds to the flow of information in the underlying problem. While here
we are not discussing any particular applications, we can still make use of that convention.
So in the case where there is a directed path from vertex $i$ to vertex $j$, where will write that
information goes from $i$ to $j$, or, more succinctly, $j$ ``sees" $i$.

The set-up of this paper is as follows. In Section \ref{chap:defns} we
give the necessary definitions concerning directed graphs, and in Section \ref{chap:Lapl}
those concerning Laplacians. In Section \ref{chap:spectra}, we discuss the spectrum of graph Laplacians,
in particular the fact that all non-zero eigenvalues have positive real part. That
means that the asymptotic behavior of the solutions of equations \eqref{eq:diffeq}
is determined by the kernel of the Laplacians.
Thus, in Section \ref{chap:kernels} we give a convenient basis for those eigenspaces. This
allows us in Section \ref{chap:dynamics} to write the asymptotics in terms of that basis.
This results in Theorem \ref{thm:Laplacian}, which is perhaps the main result in this paper.
In Section \ref{chap:confusion}, we apply this to the most important of the Laplacians,
namely the ``random walk" Laplacian. This Laplacian is particularly suited to discretization
of time, and we show that the asymptotics of the solution of the discretized equations
is again essentially the same of that of the continuous time equations. We provide examples
for all of our main statements.

Finally, a few notational issues. By ``Laplacian", we mean a matrix of the form $E-ES$,
where $E$ is diagonal with positive entries on the diagonal and $S$ is row stochastic
(details are in Section \ref{chap:Lapl}). Everything in this article goes through for
matrices of the form $E-SE$. One only needs to exchange left and right eigenvectors.
In the interest of brevity, we have not pursued this.

We use the notation ${\bf 1}_S$ for the vector whose $i$th component is 1 if $i\in S$
and 0 elsewhere. That also means that ${\bf 1}_{\{i\}}$ means the unit vector whose $i$th
component equals 1 while being 0 everywhere else. The symbol $e_i$ is used for the $i$th
diagonal element of the matrix $E$

These notes outline part of a series of 4 lectures given in summer-school/conference on mathematical
modeling of complex systems in Pescara, 2019 \cite{pescara-lectures}.
Most of this theory was described in \cite{caughveer, veerkumm} and those are the two
references that we rely most heavily on. However, various of those proofs have been substantially
simplified, and other statements have been generalized (most notably Theorem \ref{thm:leftkernel}).

\vskip 0.1in
\noindent
{\bf Acknowledgements:} JJPV is grateful to the University of Chieti-Pescara for the generous
hospitality offered.

\section{Graph Theoretic Definitions}
\label{chap:defns}

\begin{defn} A directed graph (or digraph) is a set $V=\{1,\cdots n\}$ of vertices together
with a set $E\subseteq V\times V$ or ordered pairs (the edges).
\label{def:digraph}
\end{defn}

\begin{figure}
\center
\includegraphics[width=6.0in]{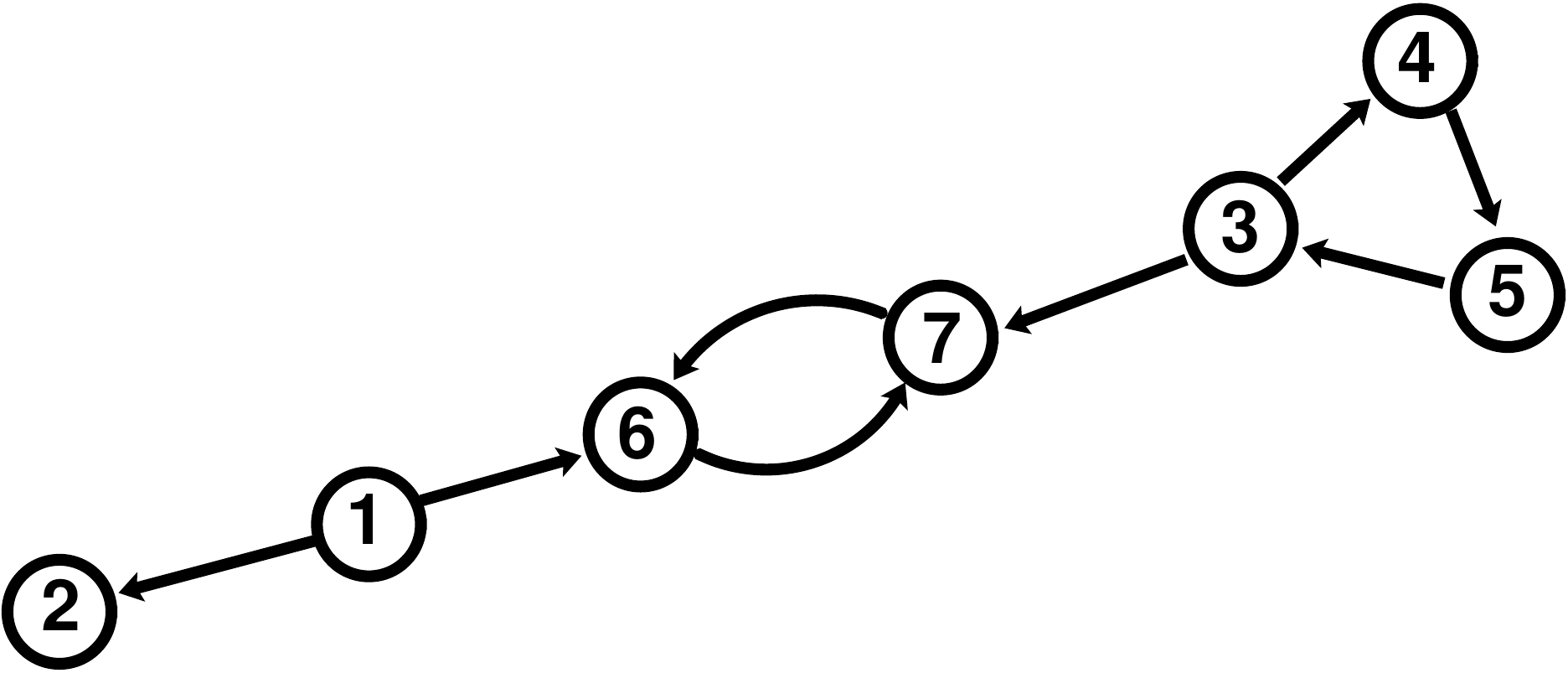}
\caption{\emph{This graph will serve as an example throughout the paper. }}
\label{fig:example}
\end{figure}

\noindent
The graph in Figure \ref{fig:example} will serve as our example of a digraph. Edges will be indicated
by $i\rightarrow j$ or $(i,j)$ (or $ij$ for short). So the graph in the figure has edges $(1,2)$, $(1,6)$,
$(6,7)$, et cetera, but it does \emph{not} have the edges $(2,1)$ and $(6,1)$.
Directed paths from $i$ to $j$ are denoted by $i\rightsquigarrow j$.
We will express this informally as: information goes from $i$ to $j$, or: $j$ ``sees" $i$.
For example, the graph in Figure \ref{fig:example} has a path $4\rightsquigarrow 6$, but there is no
path $6\rightsquigarrow 4$.

Connectedness for \emph{undirected} graph is straightforward: an undirected graph is either connected
of it is not. However, for a digraph, the situation is slightly more complicated. We need the notion
of \emph{underlying graph}. This is the undirected graph one obtains by erasing the direction
of the edge. Equivalently, it is the graph obtained by adding to each directed edge an edge in the
opposite direction.

\begin{defn} i) A digraph $G$ is strongly connected if for every ordered pair of vertices $(i,j)$,
there is a path $i\rightsquigarrow j$. Equivalently, if for every pair $i$ and $j$:
$i\leftrightsquigarrow j$.\\
 ii) A digraph $G$ is unilaterally connected if for every ordered pair of vertices $(i,j)$,
there is a path $i\rightsquigarrow j$ or a path $j\rightsquigarrow i$. \\
 iii) A digraph $G$ is weakly connected if the underlying undirected
graph is connected.\\
 iv) A digraph $G$ is not connected if it is not weakly connected.
\label{def:connected}
\end{defn}

\noindent
A subgraph which is strongly connected is called a \emph{strongly connected component}.
We will frequently abbreviate this to SCC.

The study of a graph that is not connected is of course equivalent to the study of
each its components. So the most general graph we want to study is weakly connected.
The graph of Figure \ref{fig:example} is an example of such a graph. We will need some terminology
to indicate certain subgraphs. We borrow our terminology from \cite{caughveer} and \cite{veerkumm}.

\begin{defn} i) Let $i \in V$. The \emph{reachable set} $R(i)$ consists of all $j\in V$ with $i\rightsquigarrow j$.\\
 ii) A \emph{reach} $R$ is a maximal reachable set, or a maximal unilaterally connected set.\\
 iii) A \emph{cabal} $B\subseteq R$ is the set of vertices from which the entire reach $R$ is reachable. If it is a single vertex, it is usually called a \emph{leader} or a \emph{root}.
 iv) The \emph{exclusive part} $H\subseteq R$ are those vertices in $R$ that do not ``see" vertices from
other reaches. \\
  v) The \emph{common part} $C\subseteq R$ are those vertices in $R$ that also ``see" vertices from other reaches.\\
\label{def:reaches}
\end{defn}

Note that every reach has a single non-empty cabal.
We illustrate these ideas using the graph in Figure \ref{fig:example}. That graph has two reaches,
$R_1=\{1,2,6,7\}$ and $R_2=\{3,4,5,6,7\}$. Their exclusive parts are $H_1=\{1,2\}$ and
$H_2=\{3,4,5\}$. The common parts are $C_1=C_2=\{6,7\}$. Finally, the cabals are $B_1=\{1\}$
and $B_2=\{3,4,5\}$. It is an interesting exercise to reverse the orientation of the edges
and do the taxonomy again. It is easy to see that the graph has again two reaches. But in general
the number of reaches need not be constant under orientation. Consider for example the graph
$1\rightarrow 2 \leftarrow 3$.

The relation $i \leftrightsquigarrow j$ between vertices of $G$ that defines an SCC (see Definition
\ref{def:connected}) an equivalence relation. Thus it gives a unique partition of the vertices of $G$.

\begin{figure}[pbth]
\center
\includegraphics[width=4.0in]{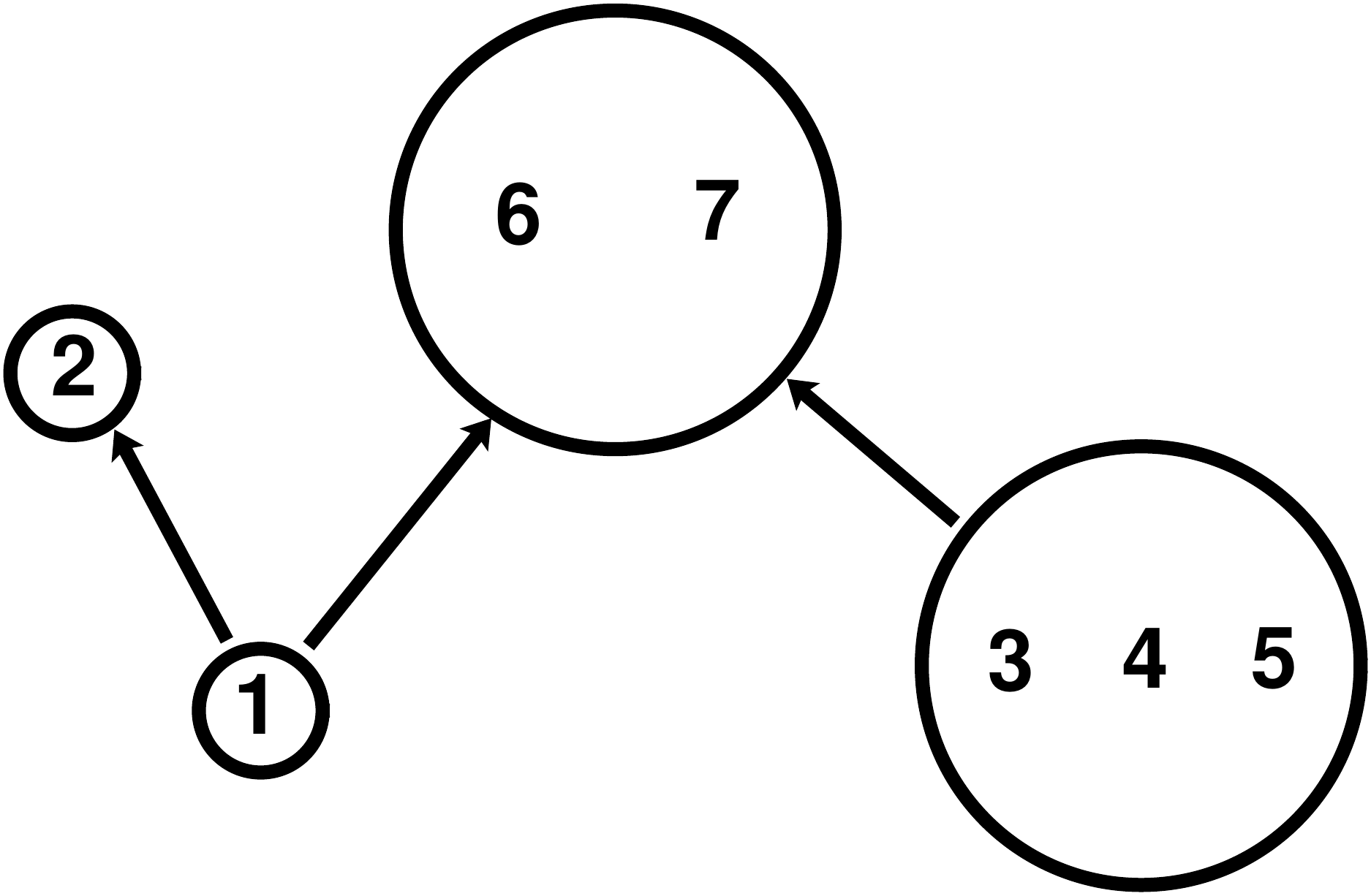}
\caption{\emph{The condensation SC$[G]$ of $G$ in Figure \ref{fig:example}.}}
\label{fig:condensation}
\end{figure}

\begin{defn} The condensation SC$[G]$ of $G$ is the graph obtained by identifying
vertices of the same SCC (or grouping them together). See Figure \ref{fig:condensation}.
\label{def:condensation}
\end{defn}

These equivalence classes respect the categories of Definition \ref{def:reaches}. For example,
given that $i$ is in a cabal, then $i\leftrightsquigarrow j$ is equivalent to $j$ is in the same
cabal. We leave it to the reader to check the other categories.
Notice that SC$[G]$ can have no cycles and therefore all the cabals are singletons.

\begin{defn} Given a digraph $G$ with vertex $i$, then ${\cal N}_i$ stands for the
set of vertices $j$ such that there is an edge $j\rightarrow i$. This is also called the (in-degree) \emph{neighborhood} of $i$.
\label{def:nbhd}
\end{defn}

\section{Graph Laplacians}
\label{chap:Lapl}

\begin{defn} The combinatorial adjacency matrix $Q$ of the
graph $G$ is defined as $Q_{ij}>0$ if there is an edge $ji$ (if ``$i$ sees $j$") and 0 otherwise.
If vertex $i$ has no incoming edges, set $Q_{ii}=1$ (create a loop).
\label{def:adjmatrix}
\end{defn}

\noindent
The last convention, on loops, is only adopted to ensure that the degree matrix $D$,
defined below, can be taken to be non-singular (and thus invertible).
One can drop the convention, but then one has to define the so-called pseudo-inverse of $D$.
This is the approach taken in \cite{caughveer}. The two approaches are equivalent.

The non-zero values of $Q_{ij}$ are the \emph{weights} of the edges $(j,i)$.
In the interest of brevity, our main example in Figure \ref{fig:example} has unit weights.
However, everything goes through in the general case.

\noindent
\begin{defn} The in-degree matrix $D$ is a diagonal matrix whose diagonal entry corresponding to
the vertex $i$ equals the sum of the weights of the edges $ji$ arriving at $i$: $d_i\equiv
\sum_jQ_{ij}$.
\label{def:indegree-matr}
\end{defn}

\noindent
The matrices $D$ and $Q$ are used to generate $S$, a row stochastic (non-negative, every row adds to 1)
version of the adjacency matrix.

\noindent
\begin{defn} The (row) stochastic matrix $S\equiv D^{-1}Q$ is called the normalized adjacency matrix.
\label{def:norm-adjmatrix}
\end{defn}

\noindent
\begin{defn} Let $E$ be a non-negative diagonal matrix. A Laplacian is a matrix of the form
$E-ES$. Common examples are: the combinatorial (comb) Laplacian,
$L^c\equiv D-DS=D-Q$, and the random walk (rw) Laplacian, ${\cal L}\equiv I-S$.
\label{def:Lapl}
\end{defn}

\noindent
Here $Q$, $D$, and $S$ are as defined earlier.

If $L$ is a Laplacian matrix, then $(Lx)_k$ is equal to $\sum_{i\neq k} \alpha_i(x_k-x_i)$
for some combination of $\alpha_i$. Thus, Laplacians describe relative observations. This
is usually called ``decentralized". Clearly, matrices with this property must have row sum zero.
It shares this property --- and hence the name --- with the discretization of the second
derivative (or ``Laplacian") of a function $f:\R\rightarrow \R$:
\bsenn
f''(j) \approx f(j-1)-2f(j)+f(j+1)
\esenn
Notice, however, that in Definition \ref{def:Lapl}, the above expression would be the negative
of a Laplacian. This convention we use ensures that Laplacians have eigenvalues whose real part
is non-negative.

As an example, we work out the matrices corresponding to the graph $G$ of Figure \ref{fig:example} assuming all weights are 1.
\bse
Q =\left( \begin {array}{c|c|ccc|cc} 1&0&0&0&0&0&0\\ \hline 1&0&0&0&0&0&0\\
\hline 0&0&0&0&1&0&0\\ 0&0&1&0&0&0&0\\ 0&0&0&1&0&0&0\\ \hline 1&0&0&0&0&0&1\\
0&0&1&0&0&1&0\end {array} \right) \quad
D=\textrm{diag}\left( \begin {array}{c} 1\\1\\1\\1\\1\\2\\2\end {array} \right)
\label{eq:adj}
\ese
\bse
L^c \equiv D-Q=\left( \begin {array}{c|c|ccc|cc} 0&0&0&0&0&0&0\\ \hline -1&1&0&0&0&0&0\\
\hline 0&0&1&0&-1&0&0\\ 0&0&-1&1&0&0&0\\ 0&0&0&-1&1&0&0\\ \hline -1&0&0&0&0&2&-1\\
0&0&-1&0&0&-1&2\end {array} \right)
\label{eq:combLapl-ex}
\ese
The spectrum of $L^c$ is: 
\bse
\left\{0,0,1, 1,3, \dfrac 32 +i\,\dfrac{\sqrt{3}}{2}, \dfrac 32 -i\,\dfrac{\sqrt{3}}{2}\right\}.
\ese

The random walk Laplacian ${\cal L} \equiv I-D^{-1}Q$ is,
\bse
{\cal L} =\left( \begin {array}{c|c|ccc|cc} 0&0&0&0&0&0&0\\ \hline -1&1&0&0&0&0&0\\
\hline 0&0&1&0&-1&0&0\\ 0&0&-1&1&0&0&0\\ 0&0&0&-1&1&0&0\\ \hline -1/2&0&0&0&0&1&-1/2\\
0&0&-1/2&0&0&-1/2&1\end {array} \right)
\label{eq:rwLapl-ex}
\ese
The spectrum of ${\cal L}$ is given by: 
\bse
\left\{0,0,\dfrac 12, 1,\dfrac 32, \dfrac 32 +i\,\dfrac{\sqrt{3}}{2}, \dfrac 32 -i\,\dfrac{\sqrt{3}}{2}\right\}.
\ese

\begin{figure}[h]
    \centering
    \begin{minipage}{0.30\textwidth}
        \centering
        \includegraphics[scale=0.30]{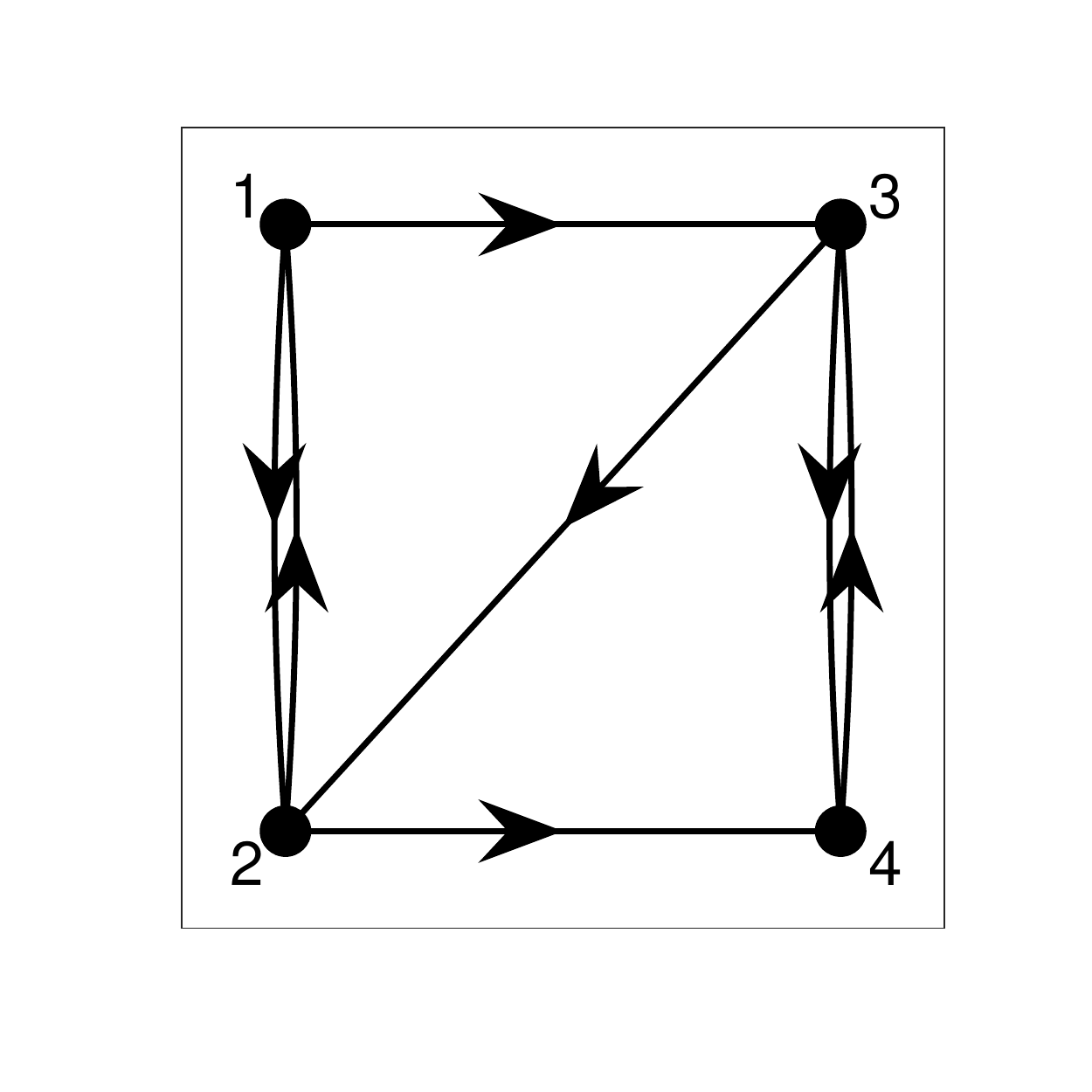}        
    \end{minipage}    
    \begin{minipage}{0.30\textwidth}
        \centering
        \includegraphics[scale=0.30]{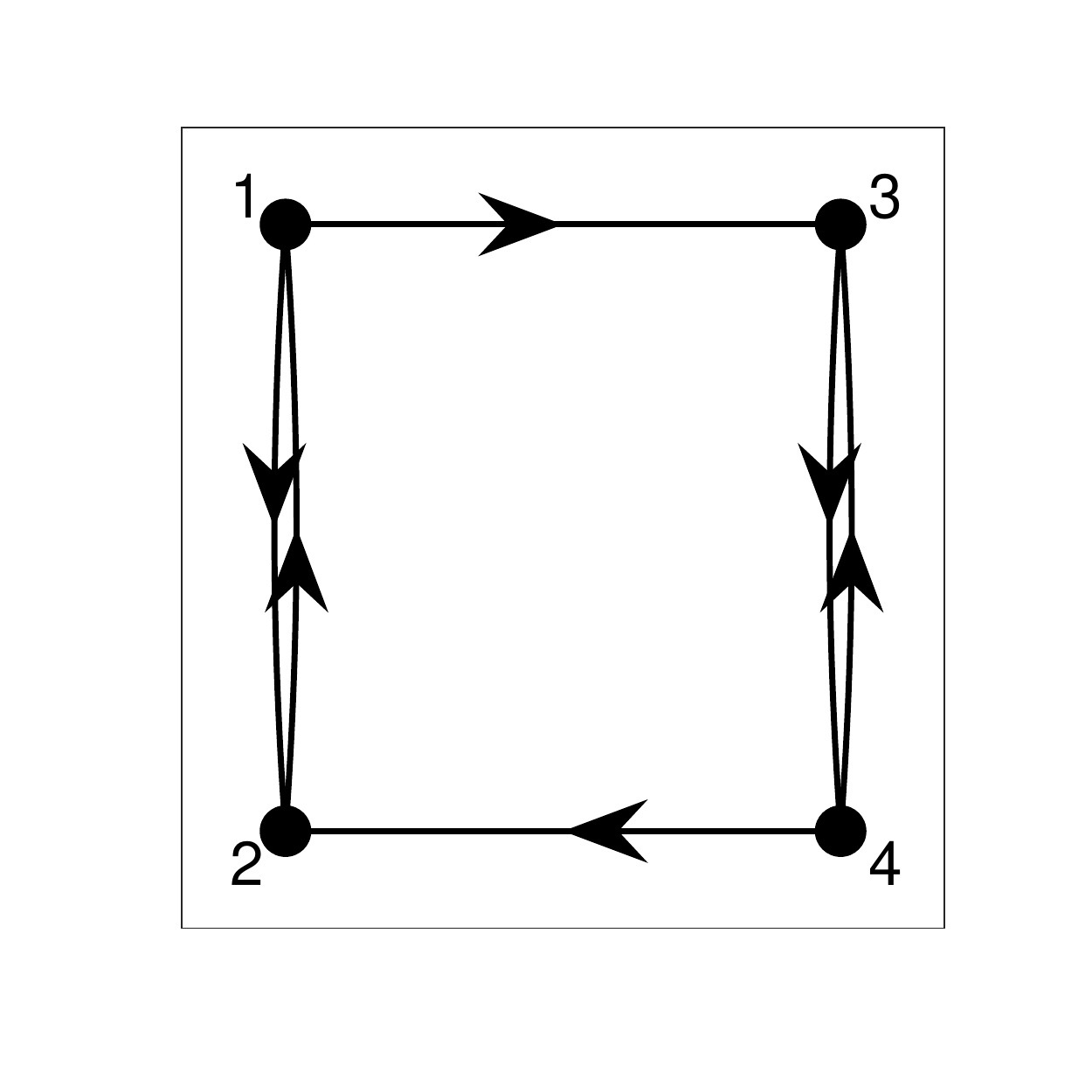}        
    \end{minipage}
    \begin{minipage}{0.30\textwidth}
       \centering
        \includegraphics[scale=0.30]{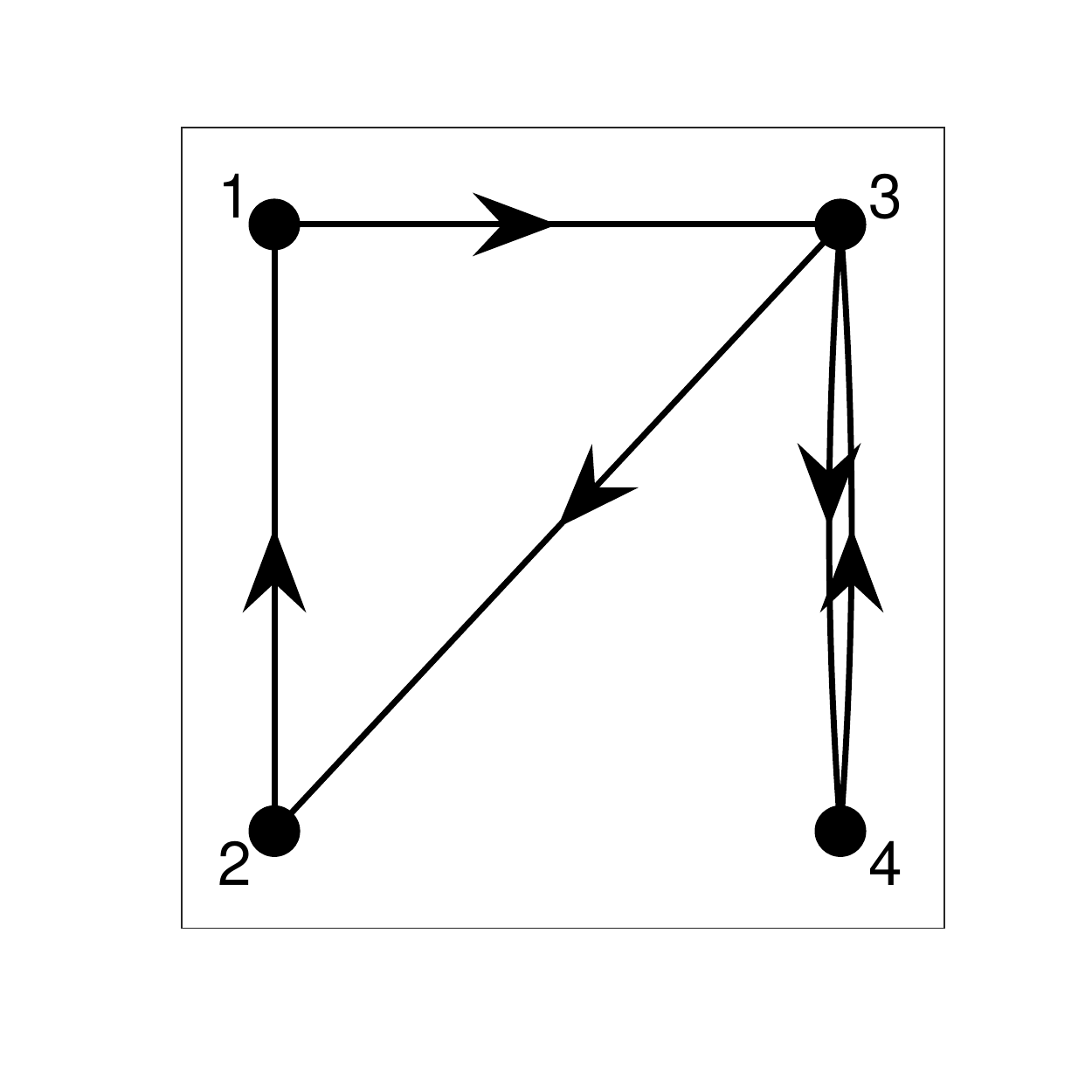}        
    \end{minipage}
    \caption{\emph{From left to right, top to bottom, three strongly connected graphs $G_1$, $G_2$, and $G_3$.}}
    \label{fig:messy}
\end{figure}

As this example shows the Laplacians do not necessarily have a real spectrum. Nor, in fact,
do they necessarily have a complete basis of eigenvectors.
We leave it to the reader to verify that in Figure \ref{fig:messy}, $L[G_1]$ and ${\cal L}[G_1]$
(with all weights equal to 1) have a non-real spectrum
and that ${\cal L}[G_2]$ and $L[G_3]$ have a non-trivial Jordan block of dimension 2.

\noindent
\begin{defn} Given a graph $G$, with $S = D^{-1} Q$.
Let $E^+$ and $E$ be non-negative diagonal matrices such that $E^+\geq E$ (entry-wise).
A generalized Laplacian is a matrix of the form $M=E^+-ES$. The matrix is strict generalized if $E^+\neq E$.
Common examples are: $M^c = E^{+}-DS$ with $E=D$ (comb), and ${\cal M} = I-ES$ with $E^+=I$ (rw).
\label{def:general-Lapl}
\end{defn}

\noindent
For us, the importance of this definition lies in the fact that the characteristic polynomial of the
Laplacian of a digraph is a product of characteristic polynomials of \emph{generalized} Laplacians
(see Proposition \ref{prop:block-triangular}). For example, in equation \eqref{eq:combLapl-ex},
the diagonal blocks are generalized comb Laplacians, and in \eqref{eq:rwLapl-ex},
they are generalized rw Laplacians.

\section{Spectra of Graph Laplacians}
\label{chap:spectra}

\begin{lem} Let $G$ be an undirected graph. The eigenvalues of a generalized Laplacian $M$ are real and
the eigenvectors form a complete basis. Neither holds necessarily even for Laplacians of a strongly connected digraph.
\label{lem:real-evals}
\end{lem}

\pf A matrix that is conjugate to a real symmetric matrix has real eigenvalues and its eigenvectors
form a complete basis. Now, recalling that $S=D^{-1}Q$ and that $Q$ is symmetric because $G$ is undirected,
we set $H:=ED^{-1}$, a diagonal matrix with \emph{positive} diagonal, and derive
\begin{align*}
 M &= E^+-ES= E^+-HQ
    \\
    &= E^+ - H^{1/2}H^{1/2}QH^{1/2}H^{-1/2}
    \\
     &=H^{1/2}\left(E^+ - H^{1/2}QH^{1/2}\right)H^{-1/2} .
\end{align*}
The last equality holds, because diagonal matrices commute. Since $E^+ - H^{1/2}QH^{1/2}$ is symmetric,
$M$ is conjugate to a symmetric matrix.

The counter-examples of Section \ref{chap:Lapl}, Figure \ref{fig:messy}, establish the second part of the this lemma.
\QED

\begin{prop} Every non-zero eigenvalue of a generalized Laplacian has positive real part.
\label{prop:gers}
\end{prop}

\vskip .1in \pf
Denote the diagonal elements of $E^+$ and $E$ by $e^+_i$ and $e_i$, respectively.
We have $e_i^+\geq e_i>0$. Apply Gersgorin's theorem \cite{HJ} to
\bsenn
(E^+-ES)v=\lambda v .
\esenn
It follows that all eigenvalues are in the union of the closed balls
\bsenn
B_{e_i}(e^+_i)= \{x\in \C\,|\, |x-e^+_i|\leq e_i\} .
\esenn
The statement follows. \QED

\vskip 0.0in
\begin{prop} The adjacency matrix of SC$[G]$ is lower block triangular after a reordering of the vertices of $G$.
\label{prop:block-triangular}
\end{prop}

\vskip 0.0in
\pf SC$[G]$ (see Definition \ref{def:condensation}) cannot contain any cycles because
the SCC's represented by vertices in a cycle of SC$[G]$ would in fact form a larger SCC.
The graph associated with SC$[G]$ can be drawn with the arrows pointing upward (see Figure \ref{fig:condensation}).
Then its vertices can be relabeled so that the vertex at the upper end (head) of an edge
is greater than the vertex at its tail. This is equivalent to saying that SC$[G]$ is lower block triangular.
\QED

\begin{cory} Any any generalized Laplacian $M=E^{+}-ES$ is lower block triangular after a reordering of the vertices of $G$.
The characteristic polynomial of $M$ is the product of the characteristic polynomials of the
$M[K_i]$ where the $K_i$ are the SCC's of $G$.
\label{cor:block-diagonal}
\end{cory}

\vskip 0.0in
\begin{prop} Let $G$ be an SCC. Any strict generalized Laplacian $M$ is non-singular.
Any Laplacian $L$ has eigenvalue 0 with geometric and algebraic multiplicity 1.
\label{prop:multiplicity}
\end{prop}

\vskip 0.0in
\pf Let $M=E^+-ES$ and suppose it has an eigenpair $\{0,v\}$. We can renormalize $v$ so that
the component with the largest modulus is $v_k=1$. $Mv=0$ gives:
\bsenn
\dfrac{e_k^+}{e_k}= \sum_{j\in {\cal N}_k}\,S_{kj}v_j ,
\esenn
where ${\cal N}_k$ stands for the \emph{neighborhood} of $k$ (see Definition \ref{def:nbhd}).
The left hand of this equality is greater than or equal to 1. The right hand is
an \emph{average} over entries with modulus less than or equal to 1.
The only way the sum can equal 1 is if $v_{j} = 1$ for all $j\in {\cal N}_k$.
Thus
\bse
e_k^+=e_k \quad \logand \quad \forall j\in {\cal N}_k: \,v_j=1 .
\label{eq:multiplicity}
\ese
Given any vertex $i$, there is a path $i\rightsquigarrow k$ ($G$ is SCC), and thus
\eqref{eq:multiplicity} holds for any vertex $i$.

The above reasoning proves that $M$ has eigenvalue 0 if and only if it is an actual Laplacian
(i.e. $E^+=E$). Furthermore, it shows that all members of the kernel of an actual Laplacian $L$ are multiples
of $\bf 1$ (the all ones vector). It remains to show that, in the case of a Laplacian $L = M$, the
algebraic multiplicity of the eigenvalues 0 equals 1.

If 0 has algebraic multiplicity $m>1$, there is a vector $x$ such that
\bsenn
L^mx=0 \quad \logand \quad L^{m-1}x\neq 0.
\esenn
This means that $L^{m-1}x = L \left( L^{m-2}x \right) = {\bf 1}$.
So $v = L^{m-2}x$ has $Lv = {\bf 1}$.
Suppose that Re$(v_i)$ is minimized at $i=k$. Then $Lv={\bf 1}$ leads to
\bsenn
e_k \textrm{Re}(v_k)=1+e_k \sum_{j}\,S_{kj} \textrm{Re}(v_j)   \ge 1 + e_{k}  \textrm{Re}(v_k),
\esenn
which is a contradiction.
\QED

\noindent
{\bf Remark:} An alternative proof is possible using the Perron Frobenius theorem \cite{boyle, Stern}
 to solve $Sv=v$.
For this, one must first note that being an SCC means that the normalized adjacency matrix $S$
is irreducible.

\begin{theo}
Given a digraph $G$. The algebraic and geometric multiplicity of the eigenvalue 0 of $L$
equals the number of reaches. A strict generalized Laplacian is non-singular. All non-zero eigenvalues have positive real part.
\label{thm:multiplicity}
\end{theo}

\vskip 0.0in
\pf We can partition the vertices the vertices of $G$ in SCC's. By Corollary \ref{cor:block-diagonal},
upon reshuffling the SCC's, the resulting Laplacian matrix is lower block triangular, and each
diagonal block is a generalized Laplacian. By Proposition \ref{prop:multiplicity},
the geometric and algebraic multiplicity of 0 equals the number of diagonal blocks (or SCC's) that
are \emph{actual} Laplacians. The generalized Laplacian of a diagonal is an actual
Laplacian if and only if that SCC has no edges coming in from other SCC's. But that happens if
and only if that SCC is a cabal. The number of cabals equals the number of reaches. \QED

\section{Kernels Right and Left}
\label{chap:kernels}

\begin{theo} Let $G$ be a digraph with $k\geq 1$ reaches.
The \emph{right} kernel of a Laplacian $L$ consists of the \emph{column} vectors
$\{\gamma_1,\cdots, \gamma_k\}$, where:
\bsenn
\left\{\begin{matrix}
\gamma_{m,j}=1 & \logif & j\in H_m &\textrm{(exclusive)}\\
\gamma_{m,j} \in (0,1) & \logif & j\in C_m &\textrm{(common)}\\
\gamma_{m,j}=0 & \logif & j\not\in R_m &\textrm{(not in reach)}\\
\sum_{m=1}^k\,\gamma_{m,j} ={\bf 1} &   &  &
\end{matrix} \right.
\esenn
\label{thm:rightkernel}
\end{theo}

\vskip .1in \pf
Pick any of the $k$ reaches and denote it by $R$. Denote its exclusive
respectively, common parts by $H$ and $C$. Recall (after Definition \ref{def:condensation})
that the SCC's respect these categories. Thus, let $X$ consist of the SCC's outside $R$ that
are ``seen" by $C$, and $Z$ the ones that are not ``seen" by the cabals in $C$.

We can chop up the Laplacian by looking at the interactions between those four groupings of SCC's:
$H$, $X$, $C$, and $Z$. For example, $X$ does not ``see" $Z$, because otherwise $C$ would also
``see" $Z$. Similarly, $X$ does not ``see", because otherwise $X$ would ``see" $H$, and therefore
be part of the same reach. This way, we obtain the schematic Laplacian given in equation
\eqref{eq:choppedlaplacian}. This matrix is block triangular in agreement with Proposition
\ref{prop:block-triangular}.

We obtain an vector in the null space of $L$ if we can solve the following equation.
\bse
\begin{pmatrix} L_{HH}  & {\bf 0} & {\bf 0} & {\bf 0} \\
                {\bf 0} & L_{XX}  & {\bf 0} & {\bf 0} \\
                L_{CH}  & L_{CX}  & L_{CC}  & {\bf 0} \\
                {\bf 0} & L_{ZX}  & {\bf 0} & L_{ZZ}  \end{pmatrix}
\begin{pmatrix}  {\bf 1}_H \\ {\bf 0}_X \\ {\bf x}_C \\ {\bf 0}_Z  \end{pmatrix}=
\begin{pmatrix}  {\bf 0}_H \\ {\bf 0}_X \\ {\bf 0}_C \\ {\bf 0}_Z  \end{pmatrix} .
\label{eq:choppedlaplacian}
\ese
But this equation boils down to
\bsenn
L_{HH}{\bf 1}_H = {\bf 0}_H \quad \logand \quad L_{CH}{\bf 1}_H+ L_{CC}{\bf x}_C = {\bf 0}_C .
\esenn
The first of these is satisfied since $L$ is a Laplacian (row-sum zero) and therefore so is $L_{HH}$.
The second of these has a unique, real solution if $L_{CC}$ is invertible. The latter is true,
because we can partition $C$ into SCC's $K_i$ in such a way that $L_{CC}$ becomes lower block triangular
(Proposition \ref{prop:block-triangular}) and the restriction of $L_{CC}$ to a block $K_i$
is a strict generalized Laplacian and so has strictly positive eigenvalues (Proposition \ref{prop:multiplicity}).  Thus $L_{CC}$ is
non-singular (Corollary \ref{cor:block-diagonal}).

Denote this real eigenvector with eigenvalue 0 by $v$. Suppose that the maximum component
is $v_n$. Then the same reasoning that leads to \eqref{eq:multiplicity} shows that if $n$
``sees" a vertex $j$, then $v_j=v_n$. Thus, since $n$ ``sees" the cabal where $v_i=1$,
we must have $v_n= 1$. Similarly, by supposing that $v_n$ is the minimum component of $v_i$,
see that $v_n= 0$. Thus all values of $v_i$ are in $[0,1]$.

Every vertex $i$ in $C$ ``sees" a vertex in $H$ (with value 1) and a vertex in
$X$ (with value 0). Thus the value at $i$ is ultimately an average collection of values
that contain both 0 and 1. Thus all entries (in $C$) are in $(0,1)$.

Finally, $\sum_{m=1}^k\,\gamma_{m,j}={\bf 1}$, because $L{\bf 1}$ must be zero, and
$\sum_{m=1}^k\,\gamma_{m,j}$ is the only combination of the $\gamma's$ that equals 1
on every vertex of the exclusive parts. \QED

As an example, we compute the basis of the null space for the Laplacian given in \eqref{eq:combLapl-ex} or \eqref{eq:rwLapl-ex} (they have the same null space).
\begin{align*}
  \gamma_1  &= \left(\begin{array}{ccccccc}1& 1& 0& 0& 0& \frac23&\frac13\end{array}\right)^T
      \quad \logand
    \\
  \gamma_2  &= \left(\begin{array}{ccccccc}0& 0 & 1 & 1 & 1& \frac13 &\frac23\end{array}\right)^T
\end{align*}

We now study the \emph{left kernel} of $L$. As a mnemonic, we use the following device:
the horizontal ``overbar" on a a vector $\bar\gamma$ indicates a (horizontal) row vector.

\noindent
\begin{theo} Let $G$ be a digraph with $k\geq 1$ reaches.
The \emph{left} kernel of Laplacian $L$ consists of the \emph{row} vectors
$\{\bar\gamma_1,\cdots, \bar\gamma_k\}$, where:
\bsenn
\left\{\begin{matrix}
\bar\gamma_{m,j}>0 & \logif & j\in B_m & \textrm{(cabal)}\\
\bar\gamma_{m,j}=0 & \logif & j\not\in B_m & \textrm{(not in cabal)}\\
\sum_{j=1}^k\,\bar\gamma_{m,j}=1 &   & & \\
\{\bar\gamma_m\}_{m=1}^k \textrm{ are orthogonal} &  & &
\end{matrix} \right.
\esenn
\label{thm:leftkernel}
\end{theo}

\vskip .0in \pf
The geometric and algebraic multiplicities of the eigenvalue 0 of $L$ equal $k$ (Theorem \ref{thm:multiplicity}). All we have to do is: find $k$ vectors $\bar\gamma_i$ in the left
kernel of $L=E(I-S)$.

For each reach $R$, we split the vertices into the cabal $B$ of $R$ and the ``rest", $X$.
We obtain an vector in the left null space of $L$ if we can solve the following equation.
\bse
\begin{pmatrix}  \bar{\bf x}_B & \bar{\bf 0}_X  \end{pmatrix}
\begin{pmatrix} L_{BB}  & {\bf 0} \\
                L_{XB} & L_{XX}  \end{pmatrix}=
\begin{pmatrix}  \bar{\bf 0}_B & \bar{\bf 0}_X  \end{pmatrix} .
\label{eq:choppedlaplacian2}
\ese
$B$ is an SCC and so by by Proposition \ref{prop:multiplicity}, this has a unique solution of the form
$\bar \gamma= (\bar \gamma_B,\bar{\bf 0}_X)$.

Set $\bar v_B=\bar \gamma_B E^{-1}$. Then $\bar v_B$ satisfies
\bsenn
\bar v_B(I_{BB}-S_{BB})=0 ,
\esenn
where $S_{BB}$ is row stochastic and irreducible. The positivity of $v_B$ follows from
Perron Frobenius. (A direct proof would take a little longer.) Thus $\bar \gamma_B$ is strictly positive.
Items ii and iv of the theorem follow after normalizing.
\QED

\vskip -0.0in
The left null spaces for the Laplacians in equations \eqref{eq:combLapl-ex} and \eqref{eq:rwLapl-ex}
are the same and are spanned by:
\bsenn
\begin{array}{ll}
  \bar\gamma_1  &=  \left(\begin{array}{ccccccc}1 & 0& 0& 0 & 0& 0& 0 \end{array}\right)
        \quad \logand
    \\
  \bar\gamma_2  &= \left(\begin{array}{ccccccc}0& 0&\frac13&\frac13&\frac13& 0& 0 \end{array}\right) .
\end{array}
\esenn

For future reference, we include this definition.

\noindent
\begin{defn} For a digraph $G$ with $n$ vertices with $k$ reaches,
we define the $n\times n$ matrix $\Gamma$ whose entries are given by:
\bsenn
\Gamma_{ij}\equiv \sum_{m=1}^k\,\gamma_{m,i} \bar\gamma_{m,j} \quad \textrm{or} \quad
\Gamma=\sum_{m=1}^k\,\gamma_{m} \otimes \bar\gamma_m
\esenn
\label{def:biggamma}
\end{defn}

\noindent
For the Laplacians given in \eqref{eq:combLapl-ex} and \eqref{eq:rwLapl-ex}, we obtain
\bse
\Gamma= \sum_{m=1}^k\,\gamma_m\otimes \bar\gamma_m=
\frac{1}{9} \, \left( \begin {array}{ccccccc} 9&0&0&0&0&0&0\\ 9&0&0&0&0&0&0\\ 0&0&3&3&3&0&0\\
0&0&3&3&3&0&0\\ 0&0&3&3&3&0&0\\ 6&0&1&1&1&0&0\\ 3&0&2&2&2&0&0 \end {array} \right) .
\label{eq:Gamma}
\ese

\noindent
\begin{lem} Let $L$ be an $n\times n$ Laplacian matrix.
The kernels of Theorems \ref{thm:rightkernel} and \ref{thm:leftkernel} can be extended to bases
of (generalized) right eigenvectors (columns), $\{\gamma_i\}_{i=1}^n$, and of (generalized) left
eigenvectors (rows), $\{\bar\gamma_i\}_{i=1}^n$, such that the matrices:
\bsenn
H = \left(\begin{array}{cccc} \gamma_1& \gamma_2 & \cdots & \gamma_n\end{array}\right)
\quad \logand \quad
\bar H = \left(\begin{array}{c} \bar\gamma_1\\ \bar\gamma_2\\ \vdots \\ \bar\gamma_n\end{array}\right)
\esenn
are inverses of one another.
\label{lem:eigenvecs}
\end{lem}

\noindent
{\bf Proof:} The right and left eigenvectors of $L$ defined in Theorems \ref{thm:rightkernel}
and \ref{thm:leftkernel} already satisfy
$\bar \gamma_i \gamma_j=\delta_{ij}$ for $i,j\in\{1,\cdots k\}$.

The extension follows directly from the Jordan Decomposition Theorem \cite{HJ}. Let $J$ be the
Jordan normal form of $L$. Then that theorem tells us that there is an invertible matrix $H$
such that $LH=H J$ or $H^{-1} L= J H^{-1}$. Right multiply the first equation
by the standard column basis vector ${\bf 1}_{\{i\}}$ to show that the $i$th column of $H$
is a generalized right eigenvector. Left multiply by ${\bf 1}_{\{i\}}^T$ to see that the $i$th row
of $H^{-1}$ is a generalized left eigenvector. \QED

\begin{defn} Let $\{\gamma_i\}_{i=1}^k$ be the (column) vectors of Theorem
\ref{thm:rightkernel} and $\{\bar\gamma_i\}_{i=1}^k$ the (row) vectors of Theorem \ref{thm:leftkernel}.
Define
\bsenn
H^0 = \left(\begin{array}{cccc} \gamma_1& \gamma_2 & \cdots & \gamma_k\end{array}\right)
\quad \logand \quad
\bar{H}^0 = \left(\begin{array}{c} \bar\gamma_1\\ \bar\gamma_2\\ \vdots \\ \bar\gamma_k\end{array}\right)\;.
\esenn
\label{defn:H0}
\end{defn}

\noindent
From Definition \ref{def:biggamma}, we now easily compute the following.

\begin{lem} $H^0\bar H^0=\Gamma$.
\label{lem:Gamma}
\end{lem}

Notice the difference between Lemmas \ref{lem:eigenvecs} and \ref{lem:Gamma}.
The matrices $H$ and $\bar H$ are both $n \times n$ so that $\bar H H = I$ implies $H \bar H = I$.
However, $H^0$ is a $n \times k$ matrix and $\bar H^{0}$ is $k \times n$ and there is no such
simplification.

\section{Laplacian Dynamics: Consensus and Diffusion}
\label{chap:dynamics}

Throughout this section, we will assume that the digraph $G$ has $k$ reaches and that $L$ is
a Laplacian of $G$.
In what follows $x$ will always stand for a column vector and $p$ for a row vector.
In this section we are interested in solving the first order Laplacian equations:
\bse
\begin{array}{ll}
  \dot x &=-Lx \with x(0)=x_0,
    \\
  \dot p &= -pL \with p(0)=p_0.
\end{array}
\label{eq:diffeq2}
\ese
The first of these equations is usually called \emph{consensus} and the second is its dual problem 
of \emph{diffusion}. We shall see below why these names are appropriate and how the solutions of
these two problems are related. We start by discussing the solutions to the concensus problem.

\vskip .1in
\begin{theo} $\lim_{t\rightarrow \infty}\,e^{-Lt}=\Gamma$.
\label{thm:Laplacian}
\end{theo}

\noindent
{\bf Proof:} Let $\{\gamma_i\}_{i=1}^k$ and $\{\bar\gamma_i\}_{i=1}^k$ as in Theorems
\ref{thm:rightkernel} and \ref{thm:leftkernel}, and then extends these sets to complete basis
of generalized) eigenvectors $\{\gamma_i\}_{i=1}^n$ and $\{\bar\gamma_i\}_{i=1}^n$ as in
Lemma \ref{lem:eigenvecs}.
Let $\lambda_i$ be the eigenvalue associated with the $i$th (generalized) eigenvector
$\gamma_i$ (or $\bar \gamma_i$). So an initial condition $x_0=y_0+z_0$ can be decomposed as
\bse
\begin{array}{lll}
  & x_0 =y_0+z_0
    \\
  & \where 
     & y_0 =\sum_{i=1}^k\,\alpha_i \gamma_i 
    \\
  & \logand 
    &z_0  =\sum_{i=k+1}^n\,\alpha_i \gamma_i \;.
\end{array}
\label{eq:initcondn}
\ese

From the standard theory of linear differential equations (see, for example, \cite{arnold}), one
easily derives that the general solution of the consensus problem $\dot x =-Lx$ is given by
\bse
x(t)=e^{-{\cal L} t}x_0=\sum_{i=1}^n\, \gamma_i e^{-\lambda_i t}\xi_i(t)\;,
\label{eq:soln-consensus}
\ese
where $\xi_i(t)$ are polynomials whose degrees are less than the size of the Jordan block corresponding
to $\lambda_i$. Furthermore, if the dimension of that Jordan block equals 1, then $\xi_i(t)=\alpha_i$.
By Theorem \ref{thm:rightkernel}, we have $\lambda_i=0$ for $i\in \{1,\cdots k\}$.
Also the zero eigenvalue has only trivial Jordan blocks and so
for $i\in \{1,\cdots k\}$, $\xi_i=1$ and $\beta_i=\alpha_i$.
By Theorem \ref{thm:multiplicity}, the $\lambda_i$, in terms with $i>k$,
have positive real parts, and so these terms converge to zero. Therefore,
substitute equation (\ref{eq:initcondn}) into equation \eqref{eq:soln-consensus}\ to get
\bsenn
\lim_{t\rightarrow \infty}\, x(t) =
\sum_{m=1}^k\,\alpha_m\,\gamma_m =y_0\;.
\esenn

Next, we determine the $\alpha_i$. Definition \ref{defn:H0} and Lemma \ref{lem:eigenvecs}
imply that
\begin{align*}
  H^0\bar H^0x_0  
    &=  H^0\left(\bar H^0 \sum_{i=1}^n\,\alpha_i \gamma_i\right)
      \\
    &= H^0\left(\sum_{i=1}^k\,\alpha_i {\bf 1}_{\{i\}} \right)
      =\sum_{i=1}^k\,\alpha_i \gamma_i=y_0 .
\end{align*}

Note the change in the upper limit in the middle equality.
Notice also that the vector  ${\bf 1}_{\{i\}}$ is a $k$-dimensional column vector,
as opposed to $\gamma_i$ which is $n$-dimensional.
The result follows from Lemma \ref{lem:Gamma}.
\QED

Since the the solutions of equation \eqref{eq:diffeq2} are given by $e^{-Lt}x(0)$ and
$p(0)e^{-Lt}$, we have the following corollary.

\begin{cory} The solutions of \eqref{eq:diffeq2} satisfy:
\bsenn
\lim_{t\rightarrow \infty} x(t)=\Gamma x_0 \quad \logand \quad \lim_{t\rightarrow \infty}
p(t)=p_0\Gamma
\esenn
\label{cor:asympt}
\end{cory}

As an example, let us consider the equations \eqref{eq:diffeq2} for the comb Laplacian $L^c$
of Figure \ref{fig:example} with initial conditions
$x_0$ and $p_0$ concentrated on vertex 7 only. Then, from \eqref{eq:Gamma}, we get
\bsenn
\begin{array}{ll}
  & \lim_{t\rightarrow \infty}\,x(t) = \Gamma x_0={\bf 0} \quad \logand
    \\
  & \lim_{t\rightarrow \infty}\,p(t) =p_0\Gamma=\frac19(3,0,2,2,2,0,0) .
\end{array}
\esenn

We return to the first part of equation \eqref{eq:diffeq2} to explain why this called the consensus
problem.
The row sum of the Laplacian is zero, so we have ${\cal L}{\bf 1}=0$.
Thus ${\bf 1}$ is in the right kernel of $L$. If the eigenvalue 0 is non-degenerate,
then from Corollary \ref{cor:asympt}, we conclude that ${\bf 1}$ is the final state and every component of the vector $x$ has the same value.
The system is, as it were, in complete agreement
or consensus. Write out the differential equation in more detail and you get
\bsenn
\dot x_i=\sum_j e_iS_{ij}(x_j-x_i) .
\esenn
Thus $\dot x_i$ is influenced by the \emph{relative} (to $x_i$ itself) positions of
$x_j$ where $ji$ is a directed edge. In terms of Definition \ref{def:nbhd}, $j$ is in the
(in-degree) neighborhood of $i$. In other words, the consensus flows in the \emph{same} direction
as the information. In our formulation, consensus describes how the information spreads
over the whole graph. Another way of saying this is that the influence of a vertex $i$ exercises
over the other vertices is described by $\Gamma {\bf 1}_{\{i\}}$.

Next we explain how do these theorems apply to the diffusion problem in the second part of equation
\eqref{eq:diffeq2}.
In this case, the vanishing of the row sum of the Laplacian
implies that $\sum_i \dot p_i=0$. Thus the sum of the components of of $p$
is preserved. Writing out the equation in full, we get
\bse
\dot p_i = \sum_j p_je_jS_{ji}-p_ie_i .
\label{eq:diffusion}
\ese
Thus, if all $p_j$ are non-negative and $p_i=0$, then $\dot p_i>0$. This means
that the the positive $p$ orthant is preserved. Since probability (or mass) is non-negative,
these two observations together mean that the second process preserves total probability
or mass. Hence the name diffusion. It is important to note that \eqref{eq:diffusion} implies
that $p_i$ is influenced by the strengths of $p_j$ where $ij$ is a directed edge.
In terms of Definition \ref{def:nbhd}, $i$ is in the (in-degree) neighborhood of $j$.
In other words, diffusion flows in the direction \emph{contrary to} the direction of the information.
Diffusion, in our formulation, tracks the source of the information. Another way of saying this is that
the influencers of a vertex $i$ are described by ${\bf 1}^T_{\{i\}}\Gamma$.

\section{The Discretization of the rw Laplacian}
\label{chap:confusion}

Perhaps the most important example of Laplacians is the rw Laplacian ${\cal L}$ of Definition
\ref{def:Lapl} (with or without weighting the edges). It is most frequently used, among other things,
to describe diffusion and consensus related problems, in discrete as well as continuous time.
In contrast with the more general Laplacian, it is particularly well-behaved if it is discretized
using time step 1. For the discrete consensus problem, we obtain
\begin{align}
  x(\ell+1)-x(\ell) &=  (S-I)x(\ell) 
       \nonumber
    \\
   & \Longrightarrow \quad x(\ell+1)=Sx(\ell) .
\label{eq:discSdyn}
\end{align}
We get a similar expression for the discrete diffusion problem,
which is usually called the \emph{random walk} problem.
Thus in this section, we will explore the following equations.
\begin{align}
  & x(\ell+1)= Sx(\ell) \with x(0)=x_0 \quad \logand 
       \nonumber
    \\
  & p(\ell+1) =p(\ell)S \with p(0)=p_0.
\label{eq:discreq2}
\end{align}

For the discrete processes, we have, as for the continuous ones, a convenient characterization
of the asymptotic behavior.
However, in the discrete case the solution $x(\ell) = S^{\ell} x_0$ might exhibit periodic behavior 
and so would not converge, as we shall see in some examples below.
Thus we study the limit of the \emph{average} instead
\bse
\frac{1}{\ell} \sum\limits_{j=0}^{\ell-1} S^{j} x_{0}, \as \ell \to \infty.
\label{eq:discravg}
\ese
If $S^{\ell} x_0$ does have a limit it will be same as the limit of the average.
Using the average, we can proceed similarly to the continuous case in section \ref{chap:dynamics}.

\vskip .1in
\begin{theo} Given a Laplacian ${\cal L}=(I-S)$ and $\Gamma$ as in Definition \ref{def:biggamma}. We have\\ $\lim_{\ell\rightarrow \infty}\,\frac 1\ell \sum_{j=0}^{\ell-1}\,S^j=\Gamma$.
\label{thm:discrLaplacian}
\end{theo}

\noindent
{\bf Proof:} The proof is very similar to that of Theorem \ref{thm:Laplacian}. The difference
is in the analogue of equation \eqref{eq:soln-consensus}.
This time, $\lambda_i$ refers to an eigenvalue of $S$, not $L$, but with otherwise the
same notation, instead of that equation, we now have
\bse
x(\ell)=S^\ell x(0)=\sum_{i=1}^n\, \gamma_i \lambda_i^\ell \xi_i(t) .
\label{eq:discrsoln}
\ese
For $i\in\{1,\cdots k\}$, we have $\lambda_i=1$ and $\xi_i(t) = \alpha_i$. The left and
right eigenspaces of the eigenvalue 1 are the same as the left and right kernels of the Laplacian
(Theorem's \ref{thm:rightkernel} and \ref{thm:leftkernel}).

By Gersgorin's theorem, the eigenvalues of $S$ are in the closed unit ball.
By Perron-Frobenius \cite{boyle, Stern}, every eigenvalue $\lambda$ not equal to 1 but with
modulus 1 has equal algebraic and geometric multiplicity. Thus if $v$ is the eigenvector
corresponding to $\lambda$,
\bsenn
\frac 1\ell \sum_{j=0}^{\ell-1}\,\lambda^\ell v = 0 .
\esenn

All other eigenvalues have modulus less than 1, and so their contribution in the sum
also vanishes.
\QED

\begin{cory} The solutions $x(\ell)$ and $p(\ell)$ of \eqref{eq:discreq2} satisfy:
\bsenn
\lim_{\ell\rightarrow \infty}\,\frac 1\ell \sum_{j=0}^{\ell-1}\,x(j)=\Gamma x_0
\quad \logand \quad \lim_{\ell\rightarrow \infty}\,\frac 1\ell \sum_{j=0}^{\ell-1}\,p(j)=p_0\Gamma .
\esenn
\label{cor:asympdisc}
\end{cory}

The fact that in the discrete case, we have to account for periodic behavior explains
why in the discrete case, we must take a limit of an \emph{average}, while in the continuous case,
it is sufficient to just take a limit (see Theorem \ref{thm:Laplacian}). Again, taking Figure
\ref{fig:example} as example with ${\cal L}$ given in equation \eqref{eq:rwLapl-ex},
consider both discrete diffusion and discrete random walk with initial condition concentrated
on vertex 7 only. As in Section \ref{chap:dynamics}, we get
\begin{align*}
  & \lim_{\ell\rightarrow \infty}\,\frac 1\ell \sum_{j=0}^{\ell-1}\,x(j)={\bf 0} \quad \logand 
    \\
  & \lim_{\ell\rightarrow \infty}\,\frac 1\ell \sum_{j=0}^{\ell-1}\,p(j)=\frac19(3,0,2,2,2,0,0) .
\end{align*}
Notice that the random walker, when it arrives at vertex 3, undergoes periodic behavior.
We effectively take the average of that behavior.

\begin{figure}[h!]
    \centering
    \begin{minipage}{0.45\textwidth}
        \centering
        \includegraphics[width=0.9\textwidth]{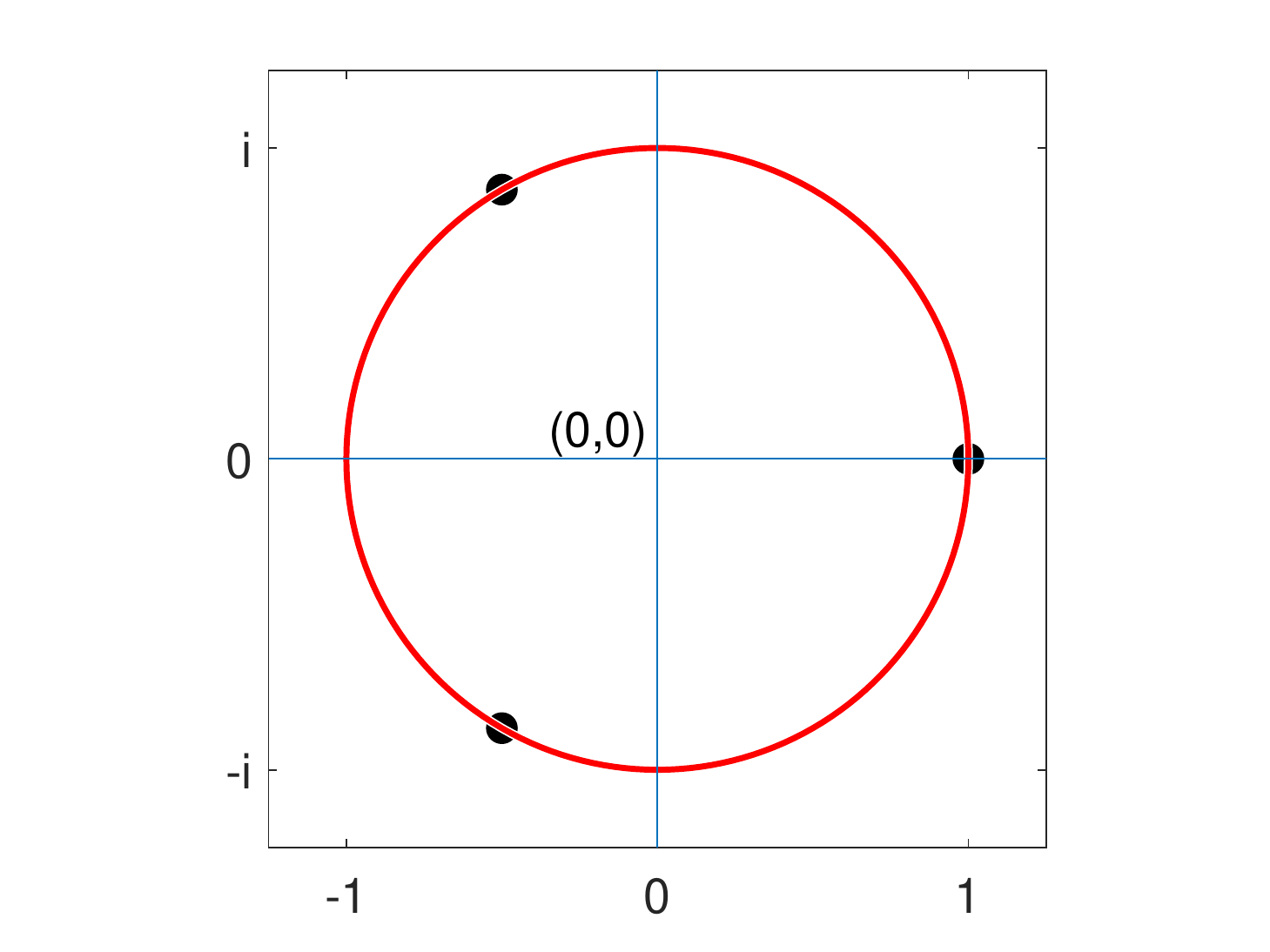} 		
        \caption{Eigenvalues of $S$ }
        \label{fig:eigvalsS}
    \end{minipage}\hfill
    \begin{minipage}{0.45\textwidth}
        \centering
        \includegraphics[width=0.9\textwidth]{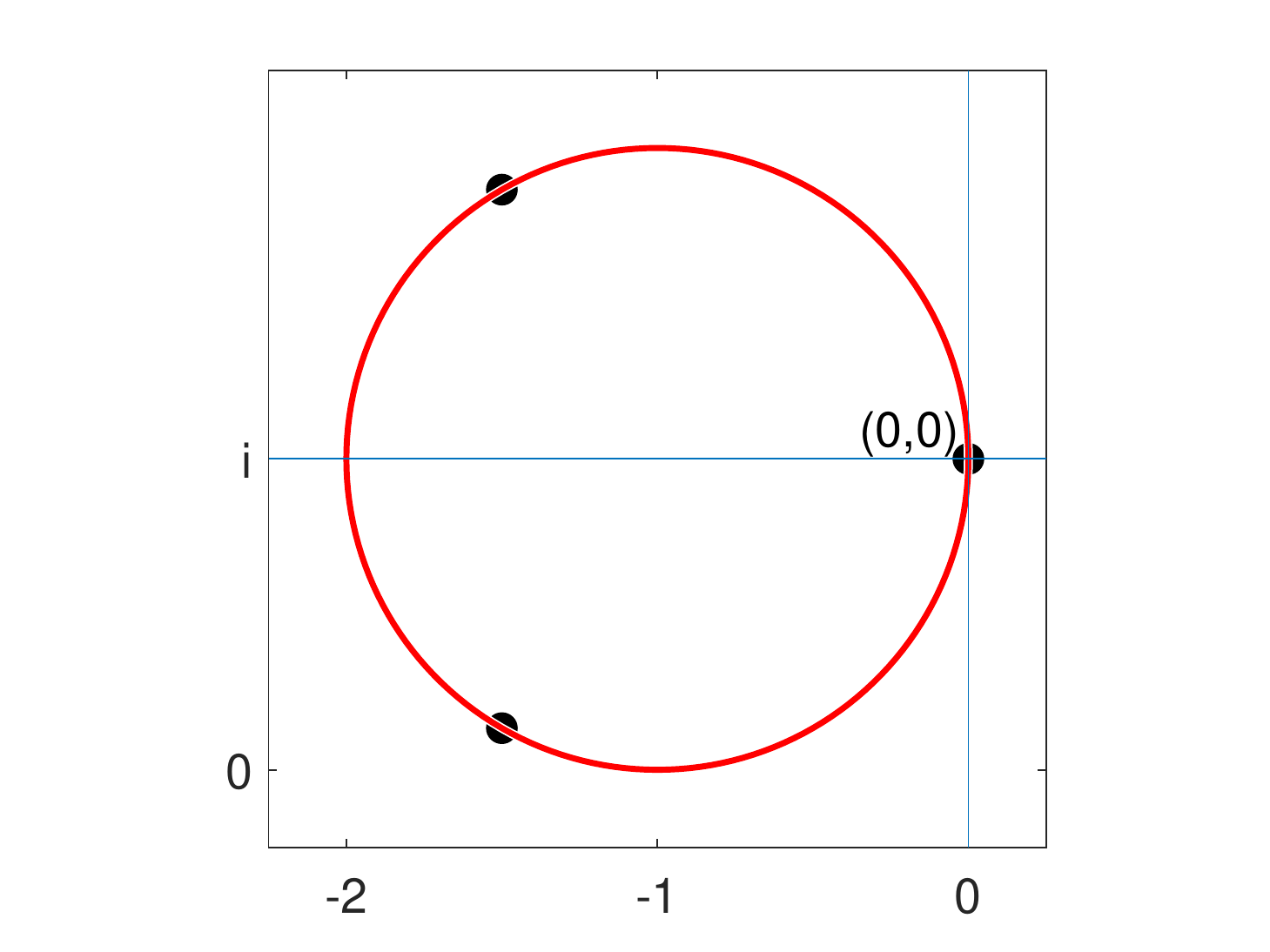}        
        \caption{Eigenvalues of $-{\cal L} = S - I$ }
        \label{fig:eigvalsL}
    \end{minipage}
\end{figure}

To isolate this periodic behavior,
take the subgraph formed by the three vertices 3, 4 and 5 of the graph in figure \ref{fig:example}.
This forms a cycle graph of order 3.
The asymptotic behavior of the discrete solutions $x( \ell ) = S^{\ell} x(0)$ are determined by the eigenvalues of $S$ that have modulus equal to $1$, because all other terms in \eqref{eq:discrsoln}
tend to 0 (since the associated eigenvalues have modulus less than 1). These are exactly
the eigenvalues of the submatrix of $S$ in \eqref{eq:discrsoln} restricted to vertices 3, 4, and 5. 
Notice that the problem of periodic behavior ``disappears" in the continuous system, because
there we consider the eigenvalues of $-{\cal L} = S - I$. So all eigenvalues shift to the left,
and all but one now have negative real part. This is illustrated if Figures \ref{fig:eigvalsS}
and \ref{fig:eigvalsL}.

Finally, we briefly discuss an alternative to the naive discretization we have been
studying in this section so far. This is the so-called time one map of the continuous
time dynamics of equation \ref{eq:diffeq2}.
\bsenn
x(\ell ) = e^{- {\cal L} } x( \ell-1 ) = \left( S^{d} \right)^{\ell } x(0),
\esenn
where we define $S^{d} = e^{- {\cal L} }$. We show that $S^{(d)}$ is a row stochastic matrix.
We start with two expansions of that matrix.
\bse
\begin{array}{lcl}
S^{d} &=& e^{-{\cal L}} = \left(I-{\cal L}+\dfrac{{\cal L}^2}{2}-\cdots\right),
    \\
S^{d} &=& e^{S-I} = e^{-1}\left(I+S+\dfrac{S^2}{2}+\cdots\right).
\end{array}
\label{eq:discretize}
\ese
From the first of these equations one easily derives that the matrix $S^{d}$ has row sum one 
(just right multiply by the vector ${\bf 1}$). The second equation implies that all its entries 
are non-negative. Therefore, $S^{d} = e^{-{\cal L}}$ is a row stochastic matrix.

It is clear that given any Laplacian ${\cal L}$ we can (in theory) always compute its
time one map $e^{- {\cal L} }$. It is interesting that the opposite is not true. From the second 
expansion in \eqref{eq:discretize}, one can deduce that for every pair of  vertices $(i,j)$
in the graph associated to the weighted adjacency matrix $S^{d}$ such that there is a path $i
\rightsquigarrow j$, there is an edge $i \rightarrow j$, though its weight might be very small. 
(A graph with this property is called \emph{transitively closed}.)
Among other things, this implies that no time one map of a Laplacian
can generate periodic behavior. In addition no Laplacian can generate
a time one map with a zero eigenvalue.

\section{Concluding Remarks}
\label{chap:conclusion}

We analyzed the asymptotic  behavior of general first order
Laplacian processes on digraphs. The most important of these are diffusion and
consensus with both continuous and discrete time. We have seen that diffusion and
consensus are \emph{dual} processes.

We remark here that given a continuous time diffusion or consensus process,
it is always possible to find its time 1 map. But vice versa is not always possible.
The reason is evident from the second part of \eqref{eq:discretize}.
That equation shows that any in time one map, every edge is realized, though not with the same weight.
More details are given in \cite{veerkumm}.

The theory presented here has more applications than anyone can write down.
We mentioned a few in the introduction. Here we want to mention briefly a few specific uses
 of the algorithms derived here. The first is that the duality described here
can be used to give a new interpretation of the famed Pagerank algorithm as described
in \cite{Stern}. The interpretation is that the Pagerank of a site corresponds
to the influence of the owner managing that site. For details see \cite{veerkumm}.

We end with two ``folklore" results that can be easily proved with the tools of this paper.
$G$ is a (weakly connected) digraph with rw Laplacian ${\cal L}$.
The union of its cabals is called $B$. Its complement is denoted as $B^c$.
First, a random walker starting at vertex $j$ has a probability $\gamma_{m,j}$ of ending up in the
$m$th cabal $B_m$. The second result is that expected hitting time $\tau(i)$ for a random walk starting
at vertex $i$ to reach (or hit) $B$, is the unique solution of
\bsenn
{\cal L}\tau={\bf 1}_{B^c} \textrm{  with  } \tau|_{B}=0 .
\esenn

\bibliographystyle{plain}
\bibliography{SurveyLaplacian}

\end{document}

%% file: inputlatex.tex
\newtheorem {theo} {\bf Theorem} [section]
\newtheorem {prop} [theo] {\bf Proposition}
\newtheorem {cory} [theo] {\bf Corollary}
\newtheorem {lem} [theo] {\bf Lemma}
\newtheorem {defn} [theo] {\bf Definition}

\newtheorem {rem} [theo] {\bf Remark}
\newcommand{\pf}{\noindent {\bf Proof. }}
\newcommand{\QED}{\hfill \CaixaPreta \vspace{6mm}}
\def\CaixaPreta{\vrule Depth0pt height6pt width6pt}

\newcommand{\qed}{\nopagebreak\hfill{\vrule width6pt height6pt depth0pt}}

\newcommand{\be}{\begin{eqnarray}}
\newcommand{\ee}{\end{eqnarray}}
\newcommand{\benn}{\begin{eqnarray*}}
\newcommand{\eenn}{\end{eqnarray*}}
\newcommand{\bse}{\begin{equation}}
\newcommand{\ese}{\end{equation}}
\newcommand{\bsenn}{\begin{displaymath}}
\newcommand{\esenn}{\end{displaymath}}
\newcommand{\logand}{\;\;{\rm and }\;\;}

\newcommand{\logif}{\;\;{\rm if }\;\;}

\newcommand{\where}{\;\;{\rm where }\;\;}
\newcommand{\with}{\;\;{\rm with }\;\;}

\newcommand{\as}{\;\;{\rm as }\;\;}

\newcommand{\C}{\mathbb{C}}

\newcommand{\R}{\mathbb{R}}


%% file: authorInfo-arXiv.tex
%

\author[1]{J. J. P. Veerman; e-mail: veerman@pdx.edu}
\author[1, 2]{R. Lyons; e-mail: rlyons@pdx.edu}

\affil[1]{\footnotesize Fariborz Maseeh Dept. of Math. and Stat., Portland State Univ.}
\affil[2]{\footnotesize Digimarc, 9405 SW Gemini Drive, Beaverton, OR, USA 97008-7192}